\documentclass[11pt]{amsart}

\usepackage{geometry}
\usepackage{todonotes}

\usepackage{amsfonts, adjustbox, amssymb, amscd}
\numberwithin{equation}{section}

\usepackage{bm}
\usepackage{verbatim}
\usepackage{mathrsfs}
\usepackage{graphicx}
\usepackage{tikz-cd}
\usepackage{subcaption}
\usepackage{listings}
\usepackage{subfiles}
\usepackage[toc,page]{appendix}
\usepackage{mathtools}
\usepackage{comment}
\usepackage{enumerate}
\usepackage{enumitem}
\usepackage[all]{xy}

\usepackage{graphicx}
\usepackage{appendix}
\usepackage{hyperref}
\hypersetup{
    colorlinks=true,
    citecolor=red,
    linkcolor=blue,
    filecolor=magenta,      
    urlcolor=red,
}
\newcommand{\bb}{\bm{b}}
\newcommand{\Mm}{{\bf{M}}}
\newcommand{\Nn}{{\bf{N}}}

\newcommand{\Dd}{{\bf{D}}}

\newcommand{\Qq}{\mathbb{Q}}

\newcommand{\Rr}{\mathbb{R}}

\newcommand{\Exc}{\operatorname{Exc}}

\newcommand{\rk}{\operatorname{rank}}

\newcommand{\Weil}{\operatorname{Weil}}

\newcommand{\Supp}{\operatorname{Supp}}

\newcommand{\Aa}{{\bf{A}}}

\newcommand{\Bb}{{\bf{B}}}
\newcommand{\Ff}{\mathcal{F}}

\newcommand{\Ii}{\Gamma}


\newcounter{parentnumber}

\newtheorem{thm}{Theorem}[section]
\newtheorem{conj}[thm]{Conjecture}
\newtheorem{cor}[thm]{Corollary}
\newtheorem{lem}[thm]{Lemma}
\newtheorem{prop}[thm]{Proposition}

\theoremstyle{definition}
\newtheorem{defn}[thm]{Definition}

\theoremstyle{definition}
\newtheorem{rem}[thm]{Remark}

\newtheorem{defthm}[thm]{Definition-Theorem}

\theoremstyle{definition}

\begin{document}

\title[On the equivalence between the effective adjunction conjectures]{On the equivalence between the effective adjunction conjectures of Prokhorov-Shokurov and of Li}
\author{Jingjun Han, Jihao Liu, and Qingyuan Xue}

\date{\today}

\begin{abstract}
Prokhorov and Shokurov introduced the famous effective adjunction conjecture, also known as the effective base-point-freeness conjecture. This conjecture asserts that the moduli component of an lc-trivial fibration is effectively base-point-free. Li proposed a variation of this conjecture, which is known as the $\Gamma$-effective adjunction conjecture, and proved that a weaker version of his conjecture is implied by the original Prokhorov-Shokurov conjecture.

In this paper, we establish the equivalence of Prokhorov-Shokurov's and Li's effective adjunction conjectures. The key to our proof is the formulation of a uniform rational polytope for canonical bundle formulas, which relies on recent developments in the minimal model program theory of algebraically integrable foliations by Ambro-Cascini-Shokurov-Spicer and Chen-Han-Liu-Xie.
\end{abstract}

\address{Department of Mathematics, Northwestern University, 2033 Sheridan Road, Evanston, IL 60208, USA}
\email{jliu@northwestern.edu}

\address{Shanghai Center for Mathematical Sciences, Fudan University, Shanghai, 200438, China}
\email{hanjingjun@fudan.edu.cn}

\address{Department of Mathematics, The University of Utah, Salt Lake City, UT 84112, USA}
\email{xue@math.utah.edu}

\maketitle

\tableofcontents

\section{Introduction}\label{sec:Introduction}

We work over the field of complex numbers $\mathbb C$.

Prokhhorov and Shokurov famously proposed the effective base-point-freeness conjecture on the moduli part of lc-trivial fibrations. 

\begin{conj}[{\cite[Conjecture 7.13]{PS09}}]\label{conj: ps09}
    Let $d$ be a positive integer and $\Ii_0$ a finite set of rational numbers. Then there exists a positive integer $I$ depending only on $d$ and $\Ii_0$ satisfying the following. Assume that
    \begin{enumerate}
        \item $f: (X,B)\rightarrow Z$ is an lc-trivial fibration such that $\dim X-\dim Z=d$, and
        \item the coefficients of the horizontal$/Z$ part of $B$ belong to $\Ii_0$.
    \end{enumerate} 
    Then $I\Mm$ is base-point-free, where $\Mm$ is the moduli part of $f: (X,B)\rightarrow Z$.
\end{conj}
Conjecture \ref{conj: ps09} is known for its complexity and has only been proven when $d=1$, as outlined in \cite[Theorem 8.1]{PS09}. The non-effective version of this conjecture for $d=2$ was recently proven in \cite[Theorem 1.4]{ABBDILW23}. However, for $d\geq 3$, Conjecture \ref{conj: ps09} remains largely unresolved.

The importance of Conjecture \ref{conj: ps09} lies in its close relationship with moduli theory. Specifically, since the moduli parts of lc-trivial fibrations characterize the moduli space of the family $X\rightarrow Z$, Conjecture \ref{conj: ps09} is crucial for the study of the moduli of varieties, especially log Calabi-Yau varieties (cf. \cite{ABBDILW23}).

Recent developments in moduli theory suggest that instead of considering only the moduli theory of pairs with standard or rational coefficients, it is more natural to include pairs with arbitrary coefficients in $[0,1]$ or $(\frac{1}{2},1]$ (cf. \cite[6.26-6.28]{Kol23}). In particular, pairs with irrational coefficients need to be considered. Since Conjecture \ref{conj: ps09} only considers lc-trivial fibrations with rational horizontal coefficients, it becomes natural to inquire whether a generalization of Conjecture \ref{conj: ps09} for lc-trivial fibrations with irrational coefficients is feasible. Fortunately, Z. Li has proposed such a variation in \cite[Conjecture 3.5(1)]{Li20}, adopting the notation of $\Ii$-base-point-freeness. In this paper, we propose a stronger version of \cite[Conjecture 3.5(1)]{Li20}.

\begin{defn}[{\cite[Definition 3.4]{Li20}}]\label{def: gamma-bpf}
    Let $\Ii\subset (0,1]$ be a set. A $\bb$-divisor $\Dd$ on a normal projective variety $X$ is called \emph{$\Ii$-base-point-free} if there exist $a_1,\dots,a_k\in\Ii$ and base-point-free $\bb$-divisors $\Dd_1,\dots,\Dd_k$, such that $\sum_{i=1}^ka_i=1$ and $\Dd=\sum_{i=1}^ka_i\Dd_i$.
\end{defn}

\begin{conj}[{cf. \cite[Conjecture 3.5(1)]{Li20}}]\label{conj: gamma-adj}
    Let $d$ be a positive integer and $\Ii\subset [0,1]$ a DCC set of real numbers. Then there exist a positive integer $I$ and a finite set $\Ii_0\subset (0,1]$ depending only on $d$ and $\Ii$ satisfying the following. Assume that
    \begin{enumerate}
        \item $f: (X,B)\rightarrow Z$ is an lc-trivial fibration such that $\dim X-\dim Z=d$, and
        \item the coefficients of the horizontal$/Z$ part of $B$ belong to $\Ii$.
    \end{enumerate} 
    Then $I\Mm$ is $\Ii_0$-base-point-free, where $\Mm$ is the moduli part of $f: (X,B)\rightarrow Z$.
\end{conj}
It is evident that Conjecture \ref{conj: gamma-adj} implies Conjecture \ref{conj: ps09}. This raises the intriguing question of whether the two conjectures are, in fact, equivalent. In support of this possibility, Z. Li introduced a less stringent form of Conjecture \ref{conj: gamma-adj} in \cite[Conjecture 3.5(2)]{Li20} and proved that Conjecture \ref{conj: ps09} implies this weaker version. However, it remains unproven whether \cite[Conjecture 3.5(2)]{Li20} implies back to Conjecture \ref{conj: ps09}. Moreover, \cite[Conjecture 3.5(2)]{Li20} is notably more complex than Conjecture \ref{conj: gamma-adj}.

In our paper, we show that Prokhorov-Shokurov's Conjecture \ref{conj: ps09} and Li's Conjecture \ref{conj: gamma-adj} are, indeed, equivalent.

\begin{thm}\label{thm: equivalence ps and gamma-adj}
  For any positive integer $d$, Conjecture \ref{conj: ps09} in relative dimension $d$ and Conjecture \ref{conj: gamma-adj} in relative dimension $d$ are equivalent. 
\end{thm}

As an immediate corollary, we prove Conjecture \ref{conj: gamma-adj} when $d=1$:

\begin{cor}\label{cor: gamma-adj dim 1}
    Conjecture \ref{conj: gamma-adj} holds when $d=1$.
\end{cor}

\noindent\textit{Idea of the proof}. The idea of the proof of Theorem \ref{thm: equivalence ps and gamma-adj} is to establish a uniform rational polytope for canonical bundle formulas (see Theorem \ref{thm: urp cbf ii} below). Roughly speaking, given an lc-trivial fibration $f: (X,B)\rightarrow Z$ with moduli part $\Mm$, we want to establish a uniform decomposition $(X,B)=\sum a_i(X,B_i)$ where $\sum a_i=1$, the horizontal$/Z$ coefficients of $B_i$ are rational, and each $f: (X,B_i)\rightarrow Z$ is an lc-trivial fibration, so that $\Mm=\sum a_i\Mm_i$. By saying ``uniform", we mean that $a_i$ and the horizontal$/Z$ coefficients of $B_i$ only depend on $\dim X-\dim Z$ and the horizontal$/Z$ coefficients of $B$. Such decomposition is already non-trivial even without the uniform condition (cf. \cite[Theorem 2.23]{JLX22}). The key ingredient we need for the proof of the existence of the \emph{uniform} decomposition is the minimal model program theory for algebraically integrable foliations, which has been established very recently \cite{ACSS21,CHLX23}.

More precisely, set $\Ff$ be the foliation induced by $f: X\rightarrow Z$ and $B^h$ the horizontal$/Z$ part of $B$. The key observation is that, if $(X,\Ff,B^h)$ is lc and $(X,B)$ satisfies Property $(*)$ (cf. \cite[Definition 2.13]{ACSS21}), then $K_{\Ff}+B^h$ is exactly the moduli part of $f: (X,B)\rightarrow Z$ \cite[Proposition 3.6]{ACSS21}. Therefore, if we can decompose $(X,\Ff,B^h)$ into Property $(*)$ foliated triples uniformly, then it automatically induces a decomposition of $\Mm$. Such decomposition is possible (Theorem \ref{thm: uniform weak acss polytope}) if we replace ``Property $(*)$" with the condition ``ACSS" (cf. \cite[Definition 7.2.3]{CHLX23}) thanks to the existence of uniform rational lc polytopes for foliations (\cite[Theorem 1.5]{DLM23}, \cite[Theorem 2.4.7]{CHLX23}). The rest of the proof is a series of changes of bases and models which preserves the moduli part of the canonical bundle formula. Despite the minimal model program applied, these arguments also rely on the fact that the moduli parts of two different canonical bundle formulas are the same provided that the corresponding lc-trivial fibrations are crepant over the generic point of the base (cf. Lemma \ref{lem: m preserved under crepant}).

\medskip

\noindent\textbf{Acknowledgement}. The authors are grateful to Zhan Li who told them Conjecture \ref{conj: gamma-adj}. The authors would also like to thank Guodu Chen, Junpeng Jiao, Fanjun Meng, and Lingyao Xie for useful discussions. The third author would like to thank his advisor Christopher D. Hacon for his constant support and many helpful discussions. The second author is affiliated with LMNS at Fudan University, and has received support from the National Key Research and Development Program of China (Grant No. 2020YFA0713200). The third author has been partially supported by NSF research grants no. DMS-1801851 and DMS-1952522, as well as a grant from the Simons Foundation (Award Number: 256202).

\section{Preliminaries}

\subsection{Notations and definitions}

We will adopt the standard notations and definitions from \cite{KM98,BCHM10} and use them freely. For foliations and generalized foliated quadruples, we follow \cite{CHLX23}, which generally aligns with the notations and definitions from \cite{CS20, ACSS21, CS21}, but there may be minor differences. For $\bb$-divisors and generalized pairs, we will follow the notations and definitions from \cite{BZ16,HL23}. For the canonical bundle formula, we will follow \cite{CHLX23}, which generally aligns with the classical definitions. For the reader's convenience, we provide the following notations and definitions that are not commonly used in the literature, or have minor differences with the classical definitions:

\begin{defn}
Let $m$ be a positive integer and $\bm{v}\in\mathbb R^m$. The \emph{rational envelope} of $\bm{v}$ is the minimal rational affine subspace of $\mathbb R^m$ which contains $\bm{v}$. For example, if $m=2$ and $\bm{v}=(\frac{\sqrt{2}}{2},1-\frac{\sqrt{2}}{2})$, then the rational envelope of $\bm{v}$ is $(x_1+x_2=1)\subset\mathbb R^2_{x_1x_2}$.
\end{defn}

\begin{defn}[Lc-trivial fibration, {\cite[Definition 11.3.1]{CHLX23}}]\label{defn: lc trivial fibration gfq}
Let $(X,B)$ be a sub-pair and $f: X\rightarrow Z$ a contraction. We say that $f: (X,B)\rightarrow Z$ is an \emph{lc-trivial fibration} if
\begin{enumerate}
\item $(X,B)$ is sub-lc over the generic point of $Z$,  
\item $K_X+B\sim_{\mathbb R,Z}0$, and
\item there exists a birational morphism $h: Y\rightarrow X$ with $K_Y+B_Y=h^*(K_X+B)$, such that $-B_Y^{\leq 0}$ is $\Rr$-Cartier and 
$$\kappa_{\sigma}(Y/Z,-B_Y^{\leq 0})=0.$$
\end{enumerate}
We remark that the classical definition of lc-trivial fibration replaces condition (3) with
\begin{enumerate}
        \item[(3')] $\rk f_*\mathcal{O}_X(\lceil\Aa^*(X,B)\rceil)=1.$
\end{enumerate} 
(cf. \cite[Condition (3) of Theorem 2]{Kaw98}, \cite[Theorem 0.2]{Amb05}, \cite[Theorem 8.3.7]{Kol07},\cite[Theorem 3.6]{FG14}). It is worth to mention that, in this paper, we shall only consider lc-trivial fibrations $f: (X,B)\rightarrow Z$ so that $B\geq 0$ over the generic point of $Z$. In this case, both (3) and (3') will automatically hold, so there will be no confusion on the notation.
\end{defn}

\begin{defn}[Discrimiant and moduli parts, cf. {\cite[Definition 2.3]{ACSS21}}]
Let $(X,B)$ be a sub-pair and $f: X\rightarrow Z$ a contraction such that $(X,B)$ is generically sub-lc$/Z$. In the following, we fix a choice of $K_X$ and a choice of $K_Z$, and suppose that for any birational morphism $g: \bar X\rightarrow X$ and $g_Z: \bar Z\rightarrow Z$, $K_{\bar X}$ and $K_{\bar Z}$ are chosen as the Weil divisors such that $g_*K_{\bar X}=K_X$ and $(g_Z)_*K_{\bar Z}=K_Z$.

Let $f': X'\rightarrow Z'$ be any contraction that is birationally equivalent to $f$ such that the induced birational maps $h: X'\dashrightarrow X$ and  $h_Z: Z'\dashrightarrow Z$ are morphisms and $Z'$ is $\Qq$-factorial. We let 
$$K_{X'}+B':=h^*(K_X+B).$$ 
For any prime divisor $D$ on $Z'$, we define
$$b_D(X',B';f):=1-\sup\left\{t\mid \left(X',B'+tf'^*D\right)\text{ is sub-lc over the generic point of } D\right\}.$$
Since being sub-lc is a property that is preserved under crepant transformations, $b_D(X,B;f)$ is independent of the choices of $X'$ and $Z'$ and is also independent of $U$.

Since $(X,B)$ is generically sub-lc$/Z$, $(X',B')$ is generically sub-lc$/Z$, so we may define
$$B_{Z'}:=\sum_{D\text{ is a prime divisor on }Z'}b_D(X,B;f)D.$$
and $$M_{X'}:=K_{X'}+B'-f'^*(K_{Z'}+B_{Z'}).$$
We call $B_{Z'}$ and $M_{X'}$ the \emph{discriminant part} and \emph{trace moduli part} of $f': (X',B')\rightarrow Z'$ respectively, and call $B_Z:=(h_Z)_*B_{Z'}$ and $M_X:=h_*B$ the \emph{discriminant part} and \emph{trace moduli part} of $f: (X,B)\rightarrow Z$ respectively.

By construction, there exist two $\bb$-divisors $\Bb$ on $Z$ and $\Mm$ on $X$, such that for any contraction $f'': X''\rightarrow Z''$ that is birationally equivalent to $f$ such that the induced birational maps $h': X''\dashrightarrow X'$ and  $h_{Z'}: Z''\dashrightarrow Z'$ are morphisms and $Z''$ is $\Qq$-factorial, $\Bb_{Z''}$ is the discriminant part of $f'': (X'',B'')\rightarrow Z''$, and $\Mm_{X''}$ is the trace moduli part of  $f'': (X'',B'')\rightarrow Z''$, where 
$$K_{X''}+B'':=h'^*(K_{X'}+B').$$ 
We call $\Mm$ the \emph{moduli part} of $f: (X,B)\rightarrow Z$ and $\Bb$ the \emph{discriminant $\bb$-divisor} of $f: (X,B)\rightarrow Z$. By construction, $\Bb$ is uniquely determined and $\Mm$ is uniquely determined for any fixed choices of $K_X$ and $K_Z$.

We say that $(X,B)$ is \emph{BP-stable} over $Z$ if $\Mm$ descends to $X$.
\end{defn}

\begin{rem}[Base moduli part]
    For any lc-trivial fibration $f: (X,B)\rightarrow Z$, the canonical bundle formula indicates that
    $$K_X+B\sim_{\mathbb R}f^*(K_Z+B_Z+\Mm^Z_Z)$$
    where $B_Z$ is the discriminant part and $\Mm^Z$ is a $\bb$-divisor. Such $\Mm^Z$ is also called the ``moduli part" of $f: (X,B)\rightarrow Z$ in many references. To avoid any confusion, we shall call such $\Mm^Z$ the \emph{base moduli part} of  $f: (X,B)\rightarrow Z$. 
    
    It is clear that, the lc-trivial fibrations, the effective base-point-freeness and the effective $\Ii$-base-point-freeness of the moduli part is equivalent to that of the base moduli part.
\end{rem}

\begin{rem}
    We can similarly define lc-trivial fibrations, discriminant part, and base moduli part for foliations. We refer the reader to \cite[Definition 11.3.1, Definition-Lemma 11.5.1]{CHLX23} for details. We do not need this in the rest part of the paper.
\end{rem}

\begin{defthm}[{\cite[Definition-Theorem 5.1.2]{CHLX23}, \cite[Definition-Theorem 6.5]{LLM23}, \cite[Theorem 2.2]{ACSS21}}]\label{defthm: weak ss reduction}
Let $X$ be a normal quasi-projective variety, $X\rightarrow U$ a projective morphism, $X\rightarrow Z$ a contraction, and $B$ an $\Rr$-divisor on $X$. Then there exist a toroidal pair  $(X',\Sigma_{X'})/U$, a log smooth pair $(Z',\Sigma_{Z'})$, and a commutative diagram
 \begin{center}$\xymatrix{
X'\ar@{->}[r]^{h}\ar@{->}[d]_{f'}& X\ar@{->}[d]^{f}\\
Z'\ar@{->}[r]^{h_Z} & Z\\
}$
\end{center}
satisfying the following.
\begin{enumerate}
\item $h$ and $h_Z$ are projective birational morphisms.
\item $f': (X',\Sigma_{X'})\rightarrow (Z',\Sigma_{Z'})$ is a toroidal contraction.
\item $\Supp(h^{-1}_*B)\cup\Supp\Exc(h)$ is contained in $\Supp\Sigma_{X'}$.
\item $X'$ has at most toric quotient singularities.
\item $f'$ is equi-dimensional.
\item $X'$ is $\Qq$-factorial klt.
\end{enumerate}
We call any such $f': (X',\Sigma_{X'})\rightarrow (Z',\Sigma_{Z'})$ (associated with $h$ and $h_Z$) which satisfies (1-6) an \emph{equi-dimensional model} of $f: (X,B)\rightarrow Z$. 
\end{defthm}

\begin{defn}[{Foliated log smooth, \cite[\S 3.2]{ACSS21}, \cite[Definition 6.2.1]{CHLX23}}]\label{defn: foliated log smooth}
Let $(X,\Ff,B)$ be a foliated triple such that $\Ff$ is algebraically integrable. We say that $(X,\Ff,B)$ is \emph{foliated log smooth} if there exists a contraction $f: X\rightarrow Z$ satisfying the following.
\begin{enumerate}
  \item $X$ has at most quotient toric singularities.
  \item $\Ff$ is induced by $f$.
  \item $(X,\Sigma_X)$ is toroidal for some reduced divisor $\Sigma_X$ such that $\Supp B\subset\Sigma_X$.  In particular, $(X,\Supp B)$ is toroidal, and $X$ is $\Qq$-factorial klt.
  \item There exists a log smooth pair $(Z,\Sigma_Z)$ such that $$f: (X,\Sigma_X)\rightarrow (Z,\Sigma_Z)$$ is an equi-dimensional toroidal contraction.
\end{enumerate}
\end{defn}

\begin{defn}[ACSS, {cf. \cite[Definitions 5.4.2, 7.2.2, 7.2.3]{CHLX23}}]\label{defn: ACSS f-triple}
Let $(X,\Ff,B)$ be a foliated triple, $G\geq 0$ a reduced divisor on $X$, and $f: X\rightarrow Z$ a projective morphism. We say that $(X,\Ff,B;G)/Z$ is \emph{weak ACSS} if
\begin{enumerate}    
\item $(X,\Ff,B)$ is lc.
\item $(Z,f(G))$ is log smooth and $G=f^{-1}(f(G))$.
\item $\Ff$ is induced by $f$ and $f$ is equi-dimensional.
\item For any reduced divisor $\Sigma\geq f(G)$ such that $(Z,\Sigma)$ is log smooth, $$(X,B+G+f^*(\Sigma-f(G)))$$ 
      is lc.
\end{enumerate}
If, in addition, 
\begin{enumerate}
\item[(5)] there exist an $\Rr$-Cartier $\Rr$-divisor $D\geq 0$ on $X$, such that  $\Supp\{B\}\subset\Supp D$, and for any reduced divisor $\Sigma\geq f(G)$ such that $(Z,\Sigma)$ is log smooth, $$(X,B+D+G+f^*(\Sigma-f(G)))$$ 
      is qdlt, and
\item[(6)] For any lc center of $(X,\Ff,B)$ with generic point $\eta$, over a neighborhood of $\eta$,
    \begin{enumerate}
      \item $\eta$ is the generic point of an lc center of $(X,\Ff,\lfloor B\rfloor)$, and
       \item $f: (X,B+G)\rightarrow (Z,f(G))$ is a toroidal morphism,
    \end{enumerate}
\end{enumerate}
then we say that $(X,\Ff,B;G)/Z$ is \emph{ACSS}. If $(X,\Ff,B;G)/Z$ is (weak) ACSS, then we say that $(X,\Ff,B)/Z$ and $(X,\Ff,B)$ are (weak) ACSS.
\end{defn}

\subsection{Results on foliations with ACSS singularities}

\begin{lem}[{\cite[Lemma 7.3.3]{CHLX23}}]\label{lem: fls imply acss}
Let $(X,\Ff,B)$ be foliated log smooth foliated triple such that $\Ff$ is algebraically integrable, $f: (X,\Sigma_X)\rightarrow (Z,\Sigma_Z)$ a contraction associated to $(X,\Ff,B)$, and $G$ the vertical$/Z$ part of $\Sigma_X$. Let $B^h$ be the horizontal$/Z$ part of $B$. Then $(X,\Ff,\Supp B^h;G)/Z$ is $\Qq$-factorial ACSS.
\end{lem}

\begin{prop}\label{prop: weak ss num 0 mmp}
    Let $(X,\Ff,B)/U$ be a $\Qq$-factorial ACSS foliated triple such that $$\kappa_{\sigma}(X/U,K_{\Ff}+B)=\kappa_{\iota}(X/U,K_{\Ff}+B)=0.$$ Then we may run a $(K_{\Ff}+B)$-MMP$/U$ with scaling of an ample$/U$ divisor, which terminates with $(X',\Ff',B')/U$ such that $K_{\Ff'}+B'\sim_{\mathbb R,U}0$.

    Moreover, for any contraction $X\rightarrow Z$ and divisor $G\geq 0$ on $X$ such that $(X,\Ff,B;G)/Z$ is ACSS, any $(K_{\Ff}+B)$-MMP$/U$ is also an MMP$/Z$, and $(X',\Ff',B';G')/Z$ is ACSS. Here $G'$ is the image of $G$ on $Z'$.
\end{prop}
\begin{proof}
    The main part of the proposition follows from \cite[Proposition 11.2.1]{CHLX23}. The moreover part follows from the cone theorem of algebraically integrable foliations (\cite[Theorem 3.9]{ACSS21}, \cite[Theorem 2.2.1]{CHLX23}) and \cite[Lemma 9.1.4]{CHLX23}. 
\end{proof}

\begin{thm}[{\cite[Theorem 1.5]{DLM23}, \cite[Theorem 2.4.7]{CHLX23}}]\label{thm: uniform rational polytope gfq}
Let $r$ be a positive integer, $v_1^0,\dots,v_m^0$ positive integers, and $\bm{v}_0:=(v_1^0,\dots,v_m^0)$. Then there exists an open set $U\ni \bm{v}_0$ of the rational envelope of $\bm{v}_0$ depending only on $r$ and $\bm{v}_0$ satisfying the following. 

Let $(X,\Ff,B=\sum_{j=1}^mv_j^0B_j)$ be an algebraically integrable foliated lc triple of rank $r$, and $B_j\geq 0$ are distinct Weil divisors. Then  $(X,\Ff,\sum_{j=1}^mv_jB_j)$ is lc for any $(v_1,\dots,v_m)\in U$.
\end{thm}

\begin{prop}[cf. {\cite[Proposition 3.6]{ACSS21}, \cite[Proposition 7.3.6]{CHLX23}}]\label{prop: weak cbf gfq}
Let $(X,\Ff,B)$ be a foliated triple, $f: X\rightarrow Z$ a contraction, and $G$ a reduced divisor on $X$ such that $(X,\Ff,B;G)/Z$ is weak ACSS. Let $\Mm$ be the moduli part of $f: (X,B)\rightarrow Z$. Then:
\begin{enumerate}
  \item $K_{\Ff}+B\sim \Mm_X$.
  \item $K_{\Ff}+B\sim_{Z}K_X+B+G.$
\end{enumerate}
\end{prop}

\begin{thm}[cf. {\cite[Theorem 11.1.5]{CHLX23}}]\label{thm: lc+weak acc=bpstable}
    Let $(X,\Ff,B)$ be an lc foliated triple, $f: X\rightarrow Z$ a contraction, and $G$ a reduced divisor on $X$ such that $(X,\Ff,B;G)/Z$ is weak ACSS. Let $\Mm$ be the moduli part of $f: (X,B+G)\rightarrow Z$. Then $\Mm$ descends to $X$.
\end{thm}
\begin{proof}
    By \cite[Theorem 11.1.5]{CHLX23}, $(X,B+G)$ is BP-stable over $Z$, so $\Mm$ descends to $X$.
\end{proof}

\begin{lem}\label{lem: m preserved under crepant}
Let $(X,B)$ and $(X',B')$ be two sub-pairs. Let $f: (X,B)\rightarrow Z$ and $f': (X',B')\rightarrow Z'$ be two lc-trivial fibrations such that $f$ and $f'$ are birationally equivalent (i.e., there exist birational maps $h: X \dashrightarrow X'$ and $h_Z: Z \dashrightarrow Z'$ such that $f' \circ h = h_Z \circ f$), and $(X,B)$ and $(X',B')$ are crepant over the generic point of $Z$. Let $\Mm,\Mm^Z$ be the moduli part and base moduli part of $f: (X,B)\rightarrow Z$ respectively, and let $\Mm',\Mm'^{Z'}$ be the moduli part and base moduli part of $f': (X',B')\rightarrow Z'$ respectively.

Then $\Mm=\Mm'$ and $\Mm^Z\sim_{\mathbb R}\Mm'^{Z'}$ for any compatible choices of $K_X$ and $K_Z$.
\end{lem}
\begin{proof}
The proof follows from the same lines of the proof of \cite[Lemma 11.4.3]{CHLX23}. 

Possibly passing to a common base and resolve indeterminacy of the induced birational map $X\dashrightarrow X'$, we may assume that $f=f'$, $X=X'$, and $Z=Z'$. Now $K_X+B=K_{X}+B'$ over the generic point of $Z$, so $B-B'$ is vertical$/Z$. Since $K_X+B\sim_{\mathbb R,Z}0$ and $K_{X}+B'\sim_{\mathbb R,Z}0$, $B-B'\sim_{\mathbb R,Z}0$, so $B-B'=f^*P$ for some $\Rr$-divisor $P$ on $Z$ (cf. \cite[Lemma 2.5]{CHL23}).

Let $B_Z$ and $B_Z'$ be the discriminant part of $f: (X,B)\rightarrow Z$ and $f: (X,B')\rightarrow Z$ respectively. By the definition of the discriminant part, $B_Z=B_Z'+P$.
Therefore,
$$\Mm_X=K_X+B-f^*(K_Z+B_Z)=K_X+B-f^*(K_Z+B_Z'+P)=K_{X}+B'-f^*(K_Z+B_Z')=\Mm'_X.$$

Since we may pass to an arbitrarily high base change, we have  $\Mm=\Mm'$. By the definition of the base moduli part, $\Mm^Z\sim_{\mathbb R}\Mm'^{Z'}$.
\end{proof}

\section{Uniform rational polytope for canonical bundle formulas}

In this section, we shall establish the existence of uniform rational polytopes for canonical bundle formulas. We will provide two versions of this uniform decomposition theorem (Theorems \ref{thm: urp cbf i} and \ref{thm: urp cbf ii}), whose statements and proofs are initially similar but different afterwards. The arguments in Theorem \ref{thm: urp cbf i} is more straightforward and clear from the point of view of uniform decomposition theorems, but we shall apply Theorem \ref{thm: urp cbf ii} to prove Theorem \ref{thm: equivalence ps and gamma-adj}.

\begin{lem}\label{lem: crepant lemma}
    Let $X$ and $X'$ be two normal quasi-projective varieties, $D=\sum_{i=1}^m v_i^0D_i$ an $\Rr$-Cartier $\Rr$-divisor on $X$ and $D'=\sum_{i=1}^m v_i^0D_i'$ an $\Rr$-Cartier $\Rr$-divisor on $X'$, such that $D$ and $D'$ are crepant, and $D_i,D_i'$ are $\Qq$-divisors. Then for any vector $\bm{v}=(v_1,\dots,v_m)$ in the rational envelope of $\bm{v}_0:=(v_1^0,\dots,v_m^0)$ in $\mathbb R^m$, $D(\bm{v}):=\sum_{i=1}^mv_iD_i$ and $D'(\bm{v}):=\sum_{i=1}^mv_iD_i'$ are crepant.
\end{lem}
\begin{proof}
We may write $D=\sum_{i=1}^cr_i\bar D_i$ and $D'=\sum_{i=1}^cr_i\bar D_i'$, such that $\bar D_i,\bar D_i'$ are $\Qq$-divisors and $r_1,\dots,r_c$ are linearly independent over $\mathbb Q$. By \cite[Lemma 5.3]{HLS19}, $\bar D_i$ and $\bar D_i'$ are $\Qq$-Cartier for each $i$. Let $p: W\rightarrow X$ and $q: W\rightarrow X'$ be a common resolution. Then 
$$\sum_{i=1}^cr_ip^*\bar D_i=p^*D=q^*D'=\sum_{i=1}^cr_iq^*\bar D_i',$$
so
$$\sum_{i=1}^cr_i(p^*\bar D_i-q^*\bar D_i')=0.$$
Thus $p^*\bar D_i=q^*\bar D_i'$ for each $i$. In particular, for any $\bm{u}=(u_1,\dots,u_c)\in\mathbb R^c$, $\bar D(\bm{u}):=\sum_{i=1}^cu_i\bar D_i$ and $\bar D'(\bm{u}):=\sum_{i=1}^cu_i\bar D_i'$ are crepant. Since for any vector $\bm{v}$ in the rational envelope of $\bm{v}_0$, there exists a unique vector $\bm{u}\in\mathbb R^c$ such that $\bar D(\bm{u})=D(\bm{v})$ and $\bar D'(\bm{u})=D'(\bm{v})$, the lemma follows.
\end{proof}

\begin{thm}\label{thm: uniform weak acss polytope}
Let $d$ be a positive integer and $v_1^0,\dots,v_m^0$ real numbers. Then there exists an open subset $U\ni \bm{v}_0$ of the rational envelope of $\bm{v}_0:=(v_1^0,\dots,v_m^0)$ in $\mathbb R^m$ depending only on $d$ and $\bm{v}_0$ satisfying the following. 

Let $(X,\Ff,B)$ be an lc foliated triple of rank $d$, $f: X\rightarrow Z$ a contraction, and $G$ a reduced divisor on $X$ such that $(X,\Ff,B;G)/Z$ is ACSS. Assume that 
\begin{itemize}
\item $B=\sum_{i=1}^mv_i^0B_i$ where $B_i\geq 0$ are Weil divisors, and
\item $B(\bm{v}):=\sum_{i=1}^mv_iB_i$ for any $\bm{v}=(v_1,\dots,v_m)\in\mathbb R^m$.
\end{itemize}
Then $(X,\Ff,B(\bm{v});G)/Z$ is weak ACSS for any $\bm{v}\in U$.
\end{thm}
\begin{proof}
We check all the conditions of Definition \ref{defn: ACSS f-triple} for $(X,\Ff,B(\bm{v});G)/Z$. Condition (1) of Definition \ref{defn: ACSS f-triple} follows from Theorem \ref{thm: uniform rational polytope gfq}. Conditions (2) and (3) of Definition \ref{defn: ACSS f-triple} are obvious. So we only need to check Condition (4) of Definition \ref{defn: ACSS f-triple}. 

We only need to show that there exists an open subset $U\ni \bm{v}_0$ of the rational envelope of $\bm{v}_0$, such that for any closed point $z\in Z$ and any reduced divisor $\Sigma\geq f(G)$ such that $(Z,\Sigma)$ is log smooth, 
$$(X,B(\bm{v})+G+f^*(\Sigma-f(G)))$$
is lc over a neighborhood of $z$ for any $\bm{v}\in U$. Possibly adding components to $\Sigma$, we may assume that $z$ is an lc center of $(Z,f(G))$. 

Let $H_1,\dots,H_{\dim Z}$ be all irreducible components of $\Sigma$. Then there exists an lc center $V_z$ of $$(X,B+G+f^*(\Sigma-f(G)))$$
whose image on $Z$ is $z$, and $V_z$ is an irreducible component of $\cap_{i=1}^{\dim Z}f^{-1}(H_i)$. Since $(X,B+G+f^*(\Sigma-f(G)))$ is qdlt, $V_z$ is normal.  Moreover, $\dim V_z=d$. We let
$$K_{V_z}+B_{V_z}(\bm{v}):=(K_X+B(\bm{v})+G+f^*(\Sigma-f(G)))|_{V_z}$$
for any $\bm{v}$ in the rational envelope of $\bm{v}_0$. By repeatedly applying inversion of adjunction, we only need to show that there exists an open subset $U$ of the rational envelope of $\bm{v}_0$, such that $(V_z,B_{V_z}(\bm{v}))$ is lc for any $\bm{v}\in U$. Let $\pi_z: W_z\rightarrow V_z$ be a $\Qq$-factorial dlt modification of $(V_z,B_{V_z}(\bm{v_0}))$ and let $$B_{W_z}(\bm{v}):=(\pi_z^{-1})_*B_{V_z}(\bm{v})+\Supp\Exc(\pi_z).$$
Then by Lemma \ref{lem: crepant lemma}, we have
$$K_{W_z}+B_{W_z}(\bm{v})=(\pi_z)^*(K_{V_z}+B_{V_z}(\bm{v}))$$
for any $\bm{v}$ in the rational envelope of $\bm{v}_0$. By repeatedly applying \cite[Corollary 3.9]{Nak16}, there exists an open subset $U$ of $\bm{v}_0$ depending only on $d$ and $\bm{v}_0$, such that $(W_z,B_{W_z}(\bm{v}))$ is lc for any $\bm{v}\in U$. Here we remark that the coefficients of $B_{W_z}(\bm{v})$ are well-controlled by \cite[Lemma 3.2]{Nak16}, so \cite[Corollary 3.9]{Nak16} can be applied. Therefore, $(V_z,B_{V_z}(\bm{v}))$ is lc for any $\bm{v}\in U$, and we are done.
\end{proof}

\begin{thm}[Uniform rational polytope for canonical bundle formula, I]\label{thm: urp cbf i}
Let $d$ be a positive integer and $v_1^0,\dots,v_m^0$ real numbers. Then there exists an open subset $U\ni \bm{v}_0$ of the rational envelope of $\bm{v}_0:=(v_1^0,\dots,v_m^0)$ in $\mathbb R^m$ depending only on $d$ and $\bm{v}_0$ satisfying the following. Assume that
\begin{itemize}
\item $f: (X,B)\rightarrow Z$ is an lc-trivial fibration such that $\dim X=d$,
\item $(X,B)$ is lc,
\item $B=\sum_{i=1}^mv_i^0B_i$ where $B_i\geq 0$ are Weil divisors,
\item $B(\bm{v}):=\sum_{i=1}^mv_iB_i$ for any $\bm{v}=(v_1,\dots,v_m)\in\mathbb R^m$,
and
\item $B_Z$ and $\Mm$ are the discriminant part and the moduli part of $f: (X,B)\rightarrow Z$ respectively.
\end{itemize}
Then:
\begin{enumerate}
\item For any $\bm{v}\in U$, $f: (X,B(\bm{v}))\rightarrow Z$ is an lc-trivial fibration with discriminant part $B_Z(\bm{v})$ and moduli part $\Mm(\bm{v})$ respectively, and $(X,B(\bm{v}))$ is lc.
\item For any vectors $\bm{v}^1,\dots,\bm{v}^k$ in $U$ and positive real numbers $a_1,\dots,a_k$ such that $\sum_{i=1}^ka_i=1$, we have
$$B_Z\left(\sum_{i=1}^ka_i\bm{v}^i\right)=\sum_{i=1}^k a_iB_Z(\bm{v}^i)\text{ and }\Mm\left(\sum_{i=1}^ka_i\bm{v}^i\right)=\sum_{i=1}^k a_i\Mm(\bm{v}^i).$$
\end{enumerate}
\end{thm}
\begin{proof}
\noindent\textbf{Step 1}. In this step, we construct an equi-dimensional model of $f: (X,B)\rightarrow Z$. We then run an MMP to achieve an auxiliary model $X''$ and construct auxiliary moduli parts $\Nn(\bm{v})$ for any $\bm{v}$ in an open neighborhood of the rational envelope of $\bm{v}_0$.

\medskip

Let $f': (X',\Sigma_{X'})\rightarrow (Z',\Sigma_{Z'})$ be an equi-dimensional model of $f: (X,B)\rightarrow Z$ associated with $h: X'\rightarrow X$ and $h_Z: Z'\rightarrow Z$. Let $\tilde{B}':=h^{-1}_*B+\Supp\Exc(h)$, $\tilde{B}'(\bm{v}):=h^{-1}_*B(\bm{v})+\Supp\Exc(h)$, and $B'^h$ and $B'^h(\bm{v})$ the horizontal$/Z'$ part of $\tilde{B}'$ and $\tilde{B}'(\bm{v})$ respectively for any $\bm{v}\in\mathbb R^m$. Let $G'$ be the vertical$/Z'$ part of $\Sigma_{X'}$.

Let $\Ff'$ be the foliation induced by $f'$. Then $(X',\Ff',B'^h)$ is foliated log smooth. By Lemma \ref{lem: fls imply acss}, $(X',\Ff',B'^h;G')/Z'$ is $\Qq$-factorial ACSS. Since
$$\kappa_{\sigma}(X'/Z,K_{\Ff'}+B'^h)=\kappa_{\sigma}(X'/Z,K_{\Ff'}+\tilde{B}')=\kappa_{\sigma}(X'/Z,K_{X'}+\tilde{B}')=\kappa_{\sigma}(X/Z,K_X+B)=0$$
and
$$\kappa_{\iota}(X'/Z,K_{\Ff'}+B'^h)=\kappa_{\iota}(X'/Z,K_{\Ff'}+\tilde{B}')=\kappa_{\iota}(X'/Z,K_{X'}+\tilde{B}')=\kappa_{\iota}(X/Z,K_X+B)=0,$$
by Proposition \ref{prop: weak ss num 0 mmp}, we may run a $(K_{\Ff'}+B'^h)$-MMP$/Z$ with scaling of an ample divisor which terminates with a model $(X'',\Ff'',B''^h)/Z$ of $(X',\Ff',B'^h)/Z$ such that $K_{\Ff''}+B''^h\sim_{\mathbb R,Z}0$, and this MMP is also an MMP$/Z'$. Denote by $f'': X''\rightarrow Z'$ the induced morphism.
Let $B''^h(\bm{v})$ and $G''$ be the strict transform of $B'^h(\bm{v})$ and $G'$ on $X''$ for any $\bm{v}\in\mathbb R^m$. By Proposition \ref{prop: weak ss num 0 mmp}, $(X'',\Ff'',B''^h;G'')/Z'$ is ACSS.

By Theorem \ref{thm: uniform rational polytope gfq}, there exists an open subset $U\ni\bm{v}_0$ in the rational envelope of $\bm{v}_0$ depending only on $d$ and $\bm{v}_0$, such that both $(X'',\Ff'',B''^h(\bm{v}))$ and $(X,B(\bm{v}))$ are lc for any $\bm{v}\in U$. By \cite[Lemma 5.3]{HLS19}, $K_X+B(\bm{v})\sim_{\mathbb R,Z}0$, so $f: (X,B(\bm{v}))\rightarrow Z$ is an lc-trivial fibration for any $\bm{v}\in U$. Since  $(X'',\Ff'',B''^h;G'')/Z'$ is ACSS, by Theorem \ref{thm: uniform weak acss polytope}, possibly shrinking $U$, we may assume that $(X'',\Ff'',B''^h(\bm{v});G'')/Z'$ is weak ACSS for any $\bm{v}\in U$. By Proposition \ref{prop: weak cbf gfq} and \cite[Lemma 5.3]{HLS19}, $$K_{X''}+B''^h(\bm{v})+G''\sim_{\mathbb R,Z'}K_{\Ff''}+B''^h(\bm{v})\sim_{\mathbb R,Z'}0$$ for any $\bm{v}\in U$, so $f'': (X'',B''^h(\bm{v})+G'')\rightarrow Z'$ is an lc-trivial fibration. Let $\Nn(\bm{v})$ be the moduli part of $f'': (X'',B''^h(\bm{v})+G'')\rightarrow Z'$. By Theorem \ref{thm: lc+weak acc=bpstable}, $\Nn(\bm{v})$ descends to $X''$ for any $\bm{v}\in U$.

\medskip

\noindent\textbf{Step 2}. In this step, we show that $\Nn(\bm{v})=\Mm(\bm{v})$ and conclude the proof of the theorem.

\medskip

 Applying the negativity lemma twice, we deduce that $K_{\Ff''}+B''^h$ and $K_X+B$ are crepant over the generic point of $Z$. By Lemma \ref{lem: crepant lemma}, $K_{\Ff''}+B''^h(\bm{v})$ and $K_X+B(\bm{v})$ are crepant over the generic point of $Z$ for any $\bm{v}\in U$. Therefore, $K_{X''}+B''^h(\bm{v})+G''$ and $K_{X}+B(\bm{v})$ are crepant over the generic point of $Z$ for any $\bm{v}\in U$. By Lemma \ref{lem: m preserved under crepant}, $\Nn(\bm{v})=\Mm(\bm{v})$ for any $\bm{v}\in U$.

By Proposition \ref{prop: weak cbf gfq}, $\Mm_{X''}= K_{\Ff''}+B''^h$ and $\Mm(\bm{v})_{X''}= K_{\Ff''}+B''^h(\bm{v})$ for any $\bm{v}\in U$. Therefore, for any vectors $\bm{v}^1,\dots,\bm{v}^k$ in $U$ and positive real numbers $a_1,\dots,a_k$ such that $\sum_{i=1}^ka_i=1$, we have
$$\Mm\left(\sum_{i=1}^ka_i\bm{v}^i\right)=\sum_{i=1}^k a_i\Mm(\bm{v}^i).$$
Let $\Mm^Z(\bm{v})$ be the base moduli part of $f: (X,B(\bm{v}))\rightarrow Z$ for any $\bm{v}\in U$. Then for any $\bm{v}\in U$, $\Mm^Z(\bm{v})$ descends to $Z'$ and
$$f''^*\Mm^Z(\bm{v})_{Z'}\sim_{\mathbb R}\Mm(\bm{v})_{X''}.$$
Therefore,
$$\Mm^Z\left(\sum_{i=1}^ka_i\bm{v}^i\right)\sim_{\mathbb R}\sum_{i=1}^k a_i\Mm^Z(\bm{v}^i)$$
and
$$\sum_{i=1}^ka_i(K_Z+B_Z(\bm{v}^i)+\Mm^Z(\bm{v}^i)_Z)\sim_{\mathbb R}K_Z+ B_Z\left(\sum_{i=1}^ka_i\bm{v}^i\right)+\Mm^Z\left(\sum_{i=1}^ka_i\bm{v}^i\right)_{Z}.$$
Thus
$$\sum_{i=1}^ka_iB_Z(\bm{v}^i)\sim_{\mathbb R}B_Z\left(\sum_{i=1}^ka_i\bm{v}^i\right).$$
By the concavity of lc thresholds,
$$\sum_{i=1}^ka_iB_Z(\bm{v}^i)\geq B_Z\left(\sum_{i=1}^ka_i\bm{v}^i\right),$$
so $$\sum_{i=1}^ka_iB_Z(\bm{v}^i)=B_Z\left(\sum_{i=1}^ka_i\bm{v}^i\right)$$
and we are done.
\end{proof}

\begin{thm}[Uniform rational polytope for canonical bundle formula, II]\label{thm: urp cbf ii}
Let $d$ be a positive integer and $v_1^0,\dots,v_m^0$ real numbers. Then there exists an open subset $U\ni \bm{v}_0$ of the rational envelope of $\bm{v}_0:=(v_1^0,\dots,v_m^0)$ in $\mathbb R^m$ depending only on $d$ and $\bm{v}_0$ 
 satisfying the following. Assume that
\begin{itemize}
\item $f: (X,B)\rightarrow Z$ is an lc-trivial fibration such that $\dim X-\dim Z=d$,
\item $B=B^h+B^v$ where $B^h$ and $B^v$ are the horizontal$/Z$ part and the vertical$/Z$ part of $B$ resepctively,
\item $B^h=\sum_{i=1}^mv_i^0B^h_i$, where $B^h_i\geq 0$ are Weil divisors,
\item $B^h(\bm{v}):=\sum_{i=1}^mv_iB^h_i$ for any $\bm{v}=(v_1,\dots,v_m)\in\mathbb R^m$,
\item $B^v=\sum_{i=1}^nu_i^0B^v_i$ for some real numbers $u_i^0$ and Weil divisors $B^v_i$,
and
\item $B_Z$ and $\Mm$ are the discriminant part and the moduli part of $f: (X,B)\rightarrow Z$ respectively.
\end{itemize}
Then there exist $\Rr$-affine functions $s_1,\dots,s_n: \mathbb R^m\rightarrow\mathbb R$ satisfying the following. Let $B^v(\bm{v}):=\sum_{i=1}^ns_i(\bm{v})B_i^v$ and $B(\bm{v}):=B^h(\bm{v})+B^v(\bm{v})$ for any $\bm{v}\in\mathbb R^m$. Then:
\begin{enumerate}
\item For any $\bm{v}\in U$, $f: (X,B(\bm{v}))\rightarrow Z$ is an lc-trivial fibration with discriminant part $B_Z(\bm{v})$ and moduli part $\Mm(\bm{v})$ respectively.
\item For any vectors $\bm{v}^1,\dots,\bm{v}^k$ in $U$ and positive real numbers $a_1,\dots,a_k$ such that $\sum_{i=1}^ka_i=1$, we have
$$B_Z\left(\sum_{i=1}^ka_i\bm{v}^i\right)=\sum_{i=1}^k a_iB_Z(\bm{v}^i)\text{ and }\Mm\left(\sum_{i=1}^ka_i\bm{v}^i\right)=\sum_{i=1}^k a_i\Mm(\bm{v}^i).$$
\end{enumerate}
\end{thm}
\begin{proof}
\noindent\textbf{Step 1}. This step is almost identical to \textbf{Step 1} of the proof of Theorem \ref{thm: urp cbf i} with minor differences. In this step, we construct an equi-dimensional model of $f: (X,B)\rightarrow Z$. We then run an MMP to achieve an auxiliary model $X''$ and construct auxiliary moduli parts $\Nn(\bm{v})$ for any $\bm{v}$ in an open neighborhood of the rational envelope of $\bm{v}_0$.

\medskip

Let $f': (X',\Sigma_{X'})\rightarrow (Z',\Sigma_{Z'})$ be an equi-dimensional model of $f: (X,B)\rightarrow Z$ associated with $h: X'\rightarrow X$ and $h_Z: Z'\rightarrow Z$. Let $\tilde{B}':=h^{-1}_*B+\Supp\Exc(h)$, $\tilde{B}'(\bm{v}):=h^{-1}_*B(\bm{v})+\Supp\Exc(h)$, and $B'^h$ and $B'^h(\bm{v})$ the horizontal$/Z'$ part of $\tilde{B}'$ and $\tilde{B}'(\bm{v})$ respectively for any $\bm{v}\in\mathbb R^m$. Let $G'$ be the vertical$/Z'$ part of $\Sigma_{X'}$.

Let $\Ff'$ be the foliation induced by $f'$. Then $(X',\Ff',B'^h)$ is foliated log smooth. By Lemma \ref{lem: fls imply acss}, $(X',\Ff',B'^h;G')/Z'$ is $\Qq$-factorial ACSS. Since
$$\kappa_{\sigma}(X'/Z,K_{\Ff'}+B'^h)=\kappa_{\sigma}(X'/Z,K_{\Ff'}+\tilde{B}')=\kappa_{\sigma}(X'/Z,K_{X'}+\tilde{B}')=\kappa_{\sigma}(X/Z,K_X+B)=0$$
and
$$\kappa_{\iota}(X'/Z,K_{\Ff'}+B'^h)=\kappa_{\iota}(X'/Z,K_{\Ff'}+\tilde{B}')=\kappa_{\iota}(X'/Z,K_{X'}+\tilde{B}')=\kappa_{\iota}(X/Z,K_X+B)=0,$$
by Proposition \ref{prop: weak ss num 0 mmp}, we may run a $(K_{\Ff'}+B'^h)$-MMP$/Z$ with scaling of an ample divisor which terminates with a model $(X'',\Ff'',B''^h)/Z$ of $(X',\Ff',B'^h)/Z$ such that $K_{\Ff''}+B''^h\sim_{\mathbb R,Z}0$, and this MMP is also an MMP$/Z'$. Denote by $f'': X''\rightarrow Z'$ the induced morphism.
Let $B''^h(\bm{v})$ and $G''$ be the strict transform of $B'^h(\bm{v})$ and $G'$ on $X''$ for any $\bm{v}\in\mathbb R^m$. By Proposition \ref{prop: weak ss num 0 mmp}, $(X'',\Ff'',B''^h;G'')/Z'$ is ACSS.

By Theorem \ref{thm: uniform rational polytope gfq}, there exists an open subset $U\ni\bm{v}_0$ in the rational envelope of $\bm{v}_0$ depending only on $d$ and $\bm{v}_0$, such that $(X'',\Ff'',B''^h(\bm{v}))$ is lc for any $\bm{v}\in U$. Since  $(X'',\Ff'',B''^h;G'')/Z'$ is ACSS, by Theorem \ref{thm: uniform weak acss polytope}, possibly shrinking $U$, we may assume that $(X'',\Ff'',B''^h(\bm{v});G'')/Z'$ is weak ACSS for any $\bm{v}\in U$. By Proposition \ref{prop: weak cbf gfq} and \cite[Lemma 5.3]{HLS19}, $$K_{X''}+B''^h(\bm{v})+G''\sim_{\mathbb R,Z'}K_{\Ff''}+B''^h(\bm{v})\sim_{\mathbb R,Z'}0$$ for any $\bm{v}\in U$, so $f'': (X'',B''^h(\bm{v})+G'')\rightarrow Z'$ is an lc-trivial fibration. Let $\Nn(\bm{v})$ be the moduli part of $f'': (X'',B''^h(\bm{v})+G'')\rightarrow Z'$. By Theorem \ref{thm: lc+weak acc=bpstable}, $\Nn(\bm{v})$ descends to $X''$ for any $\bm{v}\in U$.

\medskip

\noindent\textbf{Step 2}. In this step, we use the auxiliary moduli parts $\Nn(\bm{v})$ to construct $B(\bm{v})$, and show that $\Nn(\bm{v})$ is the moduli part of $f: (X,B(\bm{v}))\rightarrow Z$, i.e. $\Nn(\bm{v})=\Mm(\bm{v})$.

\medskip

Let $\Mm^Z(\bm{v})$ be the base moduli part of  $f'': (X'',B''^h(\bm{v})+G'')\rightarrow Z'$. Since $\Nn(\bm{v})$ descends to $X''$ for any $\bm{v}\in U$, $\Mm^Z(\bm{v})$ descends to $Z'$ for any $\bm{v}\in U$, and
$$f''^*\Mm^Z(\bm{v})_{Z'}\sim_{\mathbb R}\Nn(\bm{v})_{X''}$$ 
for any $\bm{v}\in U$. By Proposition \ref{prop: weak cbf gfq}, $\Nn_{X''}= K_{\Ff''}+B''^h$ and $\Nn(\bm{v})_{X''}= K_{\Ff''}+B''^h(\bm{v})$ for any $\bm{v}\in U$. Therefore, for any vectors $\bm{v}^1,\dots,\bm{v}^k$ in $U$ and positive real numbers $a_1,\dots,a_k$ such that $\sum_{i=1}^ka_i=1$, we have
$$\Nn\left(\sum_{i=1}^ka_i\bm{v}^i\right)=\sum_{i=1}^k a_i\Nn(\bm{v}^i).$$
Thus
$$\Mm^Z\left(\sum_{i=1}^ka_i\bm{v}^i\right)\sim_{\mathbb R}\sum_{i=1}^k a_i\Mm^Z(\bm{v}^i).$$

Let $\Mm^Z$ be the base moduli part of $f: (X,B)\rightarrow Z$. Applying the negativity lemma twice, we deduce that $K_{\Ff''}+B''^h+G''$ and $K_X+B$ are crepant over the generic point of $Z$. Thus by Lemma \ref{lem: m preserved under crepant}, $\Mm^Z$ is also the base moduli part of $f'': (X'',B''^h+G'')\rightarrow Z'$. Write $K_{X'}+B':=h^*(K_X+B)$ and let $B_{Z'}$ be the discriminant part of $f': (X',B')\rightarrow Z'$. Let $B_{Z'}''$ be the discriminant part of $f'': (X'',B''^h+G'')\rightarrow Z'$. Then
$$K_{X'}+B'\sim_{\mathbb R}f'^*(K_{Z'}+B_{Z'}+\Mm^Z_{Z'})$$
and
$$K_{X''}+B''^h+G''\sim_{\mathbb R}f''^*(K_{Z'}+B''_{Z'}+\Mm^Z_{Z'}).$$
Let $D:=B_{Z'}-B''_{Z'}$, and $p: W\rightarrow X'$ and $q: W\rightarrow X''$ a common resolution. For any $\bm{v}\in U$, let
$B'(\bm{v})$ be the unique $\Rr$-divisor on $X'$ such that
$$K_{X'}+B'(\bm{v})=p_*q^*(K_{X''}+B''^h(\bm{v})+G'')+f'^*D.$$
Since $K_{X''}+B''^h(\bm{v})+G''\sim_{\mathbb R,Z'}0,$ $K_{X'}+B'(\bm{v})-f'^*D$ is crepant to $K_{X''}+B''^h(\bm{v})+G''$. Therefore, $(K_{X'}+B'(\bm{v}))-(K_{X'}+B')$ is crepant to $B''^h(\bm{v})-B''^h$. Since $K_{\Ff''}+B''^h\sim_{\mathbb R,Z}0$, $K_{\Ff''}+B''^h(\bm{v})\sim_{\mathbb R,Z}0$ by \cite[Lemma 5.3]{HLS19}, so $B''^h(\bm{v})-B''^h\sim_{\mathbb R,Z}0$. Since $K_{X'}+B'\sim_{\mathbb R,Z}0$, $K_{X'}+B'(\bm{v})\sim_{\mathbb R,Z}0$. Let $B(\bm{v}):=h_*B'(\bm{v})$. Then $K_{X}+B(\bm{v})\sim_{\mathbb R,Z}0$ for any $\bm{v}\in U$. By Lemma \ref{lem: crepant lemma}, $(X,B(\bm{v}))$ and $(X'',B''^h(\bm{v})+G'')$ are crepant over the generic point of $Z$ for any $\bm{v}\in U$. Thus $f: (X,B(\bm{v}))\rightarrow Z$ is an lc-trivial fibration for any $\bm{v}\in U$. Moreover, by Lemma \ref{lem: m preserved under crepant}, $\Mm^Z(\bm{v})$ is the base moduli part of $f: (X,B(\bm{v}))\rightarrow Z$ and $\Nn(\bm{v})=\Mm(\bm{v})$ for any $\bm{v}\in U$.

\medskip

\noindent\textbf{Step 3}. We conclude the proof in this step. 

\medskip

By our construction, $B(\cdot): U\rightarrow\Weil_{\mathbb R}(X)$ is an affine function. Let $B_Z(\bm{v})$ be the discriminant part of $f: (X,B(\bm{v}))\rightarrow Z$ for any $\bm{v}\in U$. Since
$$K_X+B(\bm{v})\sim_{\mathbb R}f^*(K_Z+B_Z(\bm{v})+\Mm^Z(\bm{v})_Z)$$
for any $\bm{v}\in U$, for any vectors $\bm{v}^1,\dots,\bm{v}^k$ in $U$ and positive real numbers $a_1,\dots,a_k$ such that $\sum_{i=1}^ka_i=1$, we have
$$\sum_{i=1}^ka_i(K_Z+B_Z(\bm{v}^i)+\Mm^Z(\bm{v}^i)_Z)\sim_{\mathbb R}K_Z+ B_Z\left(\sum_{i=1}^ka_i\bm{v}^i\right)+\Mm^Z\left(\sum_{i=1}^ka_i\bm{v}^i\right)_{Z}.$$
Thus
$$\sum_{i=1}^ka_iB_Z(\bm{v}^i)\sim_{\mathbb R}B_Z\left(\sum_{i=1}^ka_i\bm{v}^i\right).$$
By the concavity of lc thresholds,
$$\sum_{i=1}^ka_iB_Z(\bm{v}^i)\geq B_Z\left(\sum_{i=1}^ka_i\bm{v}^i\right),$$
so $$\sum_{i=1}^ka_iB_Z(\bm{v}^i)=B_Z\left(\sum_{i=1}^ka_i\bm{v}^i\right)$$ and we are done.
\end{proof}

\begin{rem}
   The same lines of the arguments of Theorems \ref{thm: urp cbf i} and \ref{thm: urp cbf ii} also work for lc-trivial fibrations of generalized pairs as all results in \cite{CHLX23} can still be applied. (For the proof of Theorem \ref{thm: uniform weak acss polytope} which uses \cite[Corollary 3.9]{Nak16} for pairs, it can be replaced with \cite[Theorem 3.6]{Che23} for generalized pairs.) Due to technicality of the arguments, we omit the detailed statements and proofs here.
\end{rem}

\begin{rem}
    We can also deduce a uniform rational polytope for canonical bundle formulas for lc-trivial fibrations with DCC coefficients, similar to \cite[Theorem 4.1]{Li20}, \cite[Theorems 1.1, 1.4]{CHL22}, and \cite[Theorem 1.9]{CHL23}. The proof almost follows from the same lines of the proofs of Theorems \ref{thm: urp cbf i} and \ref{thm: urp cbf ii}. Again, due to technicality of the arguments, we omit the detailed statements and proofs here.
\end{rem}

\section{Proof of the main theorem}

\begin{proof}[Proof of Theorem \ref{thm: equivalence ps and gamma-adj}]
    It is obvious that Conjecture \ref{conj: gamma-adj} in relative dimension $d$ implies Conjecture \ref{conj: ps09} in relative dimension $d$, so we only need to show that Conjecture \ref{conj: ps09} in relative dimension $d$ implies  Conjecture \ref{conj: gamma-adj} in relative dimension $d$. 
    
    Conditions and notations as in Conjecture \ref{conj: gamma-adj} and assume that Conjecture \ref{conj: ps09} holds in relative dimension $d$. By \cite[Theorem 1.5]{HMX14}, we may assume that $\Ii$ is a finite set $\{v_1^0,\dots,v_m^0\}$ for some non-negative integer $m$. Let $B^h$ be the horizontal$/Z$ part of $B$ and write $B^h=\sum_{i=1}^m v_i^0B_i^h$, where $v_i^0\in\Ii$ for each $i$ and $B_i^h\geq 0$ are Weil divisors. Let $B^h(\bm{v}):=\sum_{i=1}^mv_iB_i^h$ for any $\bm{v}:=(v_1,\dots,v_m)\in\mathbb R^m$. Let $\bm{v}_0:=(v_1^0,\dots,v_m^0)$. Let $B_Z$ and $\Mm$ be the discriminant part and the moduli part of $f: (X,B)\rightarrow Z$ respectively.
    
    Let $U\ni\bm{v}_0$ be an open subset of the rational envelope of $\bm{v}_0$ as in Theorem \ref{thm: urp cbf ii} which only depends only on $d$ and $\bm{v}_0$. Let $k:=\dim U+1$ and $\bm{v}^1,\dots,\bm{v}^{k}\in U\cap\mathbb Q^m$ vectors depending only on $d$ and $\bm{v}_0$, such that $\bm{v}_0$ is contained in the interior of the convex hull spanned by $\bm{v}^1,\dots,\bm{v}^{k}$. Then there exist unique real numbers $a_1,\dots,a_k\in (0,1]$ such that $\sum_{i=1}^ka_i=1$ and $\sum_{i=1}^ka_i\bm{v}^i=\bm{v}_0$. By Theorem \ref{thm: urp cbf ii}, there exist $\Rr$-divisors $B(\bm{v}^1),\dots,B(\bm{v}^k)$, such that
    \begin{itemize}
        \item $B^h(\bm{v}^i)$ is the horizontal$/Z$ part of $B(\bm{v}^i)$ for each $i$,
        \item $f: (X,B(\bm{v}^i))\rightarrow Z$ is an lc-trivial fibration for each $i$ with discriminant part $B_Z(\bm{v})$ and moduli part $\Mm(\bm{v})$ respectively, and
\item $B_Z=\sum_{i=1}^k a_iB_Z(\bm{v}^i)$ and $\Mm=\sum_{i=1}^k a_i\Mm(\bm{v}^i).$
    \end{itemize}
By Conjecture \ref{conj: ps09} in relative dimension $d$, there exists a positive integer $I$ depending only on $d$ and $\Ii$, such that $I\Mm(\bm{v}^i)$ is base-point-free. Therefore, $I\Mm$ is $\{a_1,\dots,a_k\}$-base-point-free. Let $\Ii_0:=\{a_1,\dots,a_k\}$, and Conjecture \ref{conj: gamma-adj} in relative dimension $d$ follows.
\end{proof}

\begin{proof}[Proof of Corollary \ref{cor: gamma-adj dim 1}]
    It follows from Theorem \ref{thm: equivalence ps and gamma-adj} and \cite[Theorem 8.1]{PS09}.
\end{proof}

\end{document}